\documentclass[10pt,a4paper]{article}
\usepackage[margin=2cm]{geometry}
\usepackage{amsmath,amsthm,amsfonts,amssymb,amscd,cite,graphicx}
\usepackage{latexsym}
\usepackage{enumitem}
\usepackage[usenames, dvipsnames]{color}
\definecolor{mygray}{gray}{0.7}

\usepackage{titlesec}
\titleformat{\section}
{\normalfont\fontsize{16}{20}\bfseries}{\thesection}{1em.}{}

\allowdisplaybreaks[4]

\let\oldbibliography\thebibliography
\renewcommand{\thebibliography}[1]{%
  \oldbibliography{#1}%
  \setlength{\itemsep}{-2pt}%
}
\usepackage{tikz} 

\usetikzlibrary{decorations.pathreplacing}
\usetikzlibrary{arrows,decorations.pathmorphing,snakes}
\usetikzlibrary{fit}
\usepackage{mathrsfs}

\begin{document}

\baselineskip=0.20in

\makebox[\textwidth]{%
\hglue-15pt
\begin{minipage}{0.6cm}	
%\vfill
\vskip9pt
\end{minipage} \vspace{-\parskip}
\begin{minipage}[t]{11cm}
\footnotesize{ XXXXXXXXXXXXXXXXXXXXXX
%	{\bf Enumerative Combinatorics and Applications} \\ \underline{ecajournal.haifa.ac.il}
}
\end{minipage}

\hfill
\begin{minipage}[t]{5.2cm}
\normalsize 
%{\it ECA}  {\bf X} (20XX) XX--XX
XXXXXXXXXXX
\end{minipage}}

\vskip36pt

\begin{center}
{\large \bf Philippe Flajolet's early work in combinatorics}\\[10pt]

Helmut Prodinger$^\ast$\\[20pt]

\footnotesize {\it $^\ast$Department of Mathematical Sciences,
	Stellenbosch University,
	7602 Stellenbosch,	South Africa, and\\
	Department of Mathematics and Mathematical Statistics,
	Umea University,
	907 36 Umea, 	 Sweden  \\
Email: hproding@sun.ac.za}\\[10pt]

\noindent {\footnotesize {\bf Received}: Day Month 20XX},
{\footnotesize {\bf Revised}: Day Month 20XX} {\footnotesize {\bf Accepted}: Day Month 20XX}, {\footnotesize {\bf Published}: Day Month 20XX}\\

\noindent The authors: Released under the CC BY-ND license (International 4.0)
\end{center}

\setcounter{page}{1} \thispagestyle{empty}

\baselineskip=0.30in

\normalsize

\begin{abstract}

Some of Philippe Flajolet's combinatorial contributions that he wrote between 1976 and 1995, say, are described. In most of Flajolet's papers,  asymptotic/analytic
considerations play a major role. To be true to the spirit of the journal ECA, the emphasis is on the \emph{combinatorial} part. 

Covered are Register function of binary trees, approximate counting, digital search trees and the Rice formula, probabilistic counting, adding a new slice, Ramanujan's influence,  finite differences and harmonic numbers, Digits and the Mellin-Perron formula.

\bigskip



\end{abstract}

\section{Introduction}

Philippe Flajolet (1948--2011) was, among many other things, a very powerful combinatorialist. 
In 1998, on the occasion of his 50th birthday, a paper with the title \emph{Philippe Flajolet's research in Combinatorics and Analysis of Algorithms} was published
\cite{P-Flajolet-50}. An attempt was then made to be as complete as possible and say at least of few words about all his works. 
Here, the approach is different: a few prominent themes will be discussed in \emph{more} detail. They were all very influential and highly cited, and were written in Flajolet's younger years. The author had the privilege to follow the development of these papers when they were written, and they were all studied back in Vienna with interested students.

Flajolet, after a few years in Logic and Linguistics, moved around 1976 to the \emph{Analysis of Algorithms}; if I am not mistaken, Jean Vuillemin brought the problem about the register-function (Horton-Strahler numbers) to him. After that, Flajolet wrote a series of very powerful papers, and they proved how strong he was as a problem solver. Typically, his papers have a combinatorial first part with an asymptotic analysis part to follow. In the spirit of this journal, we concentrate on the combinatorial part.
In France, at the period, there was M. P. Sch\"utzenberger as a very influential figure, and, following him, it was not uncommon for young french researchers to deal with power series, words, trees, lattice paths, and so on. The use of power series was, however, more of a formal type, and even involved non-commuting variables.

Flajolet told me on several occasions that Comtet's book \cite{comtet-book} was one of his early inspirations, and he liked the \emph{example-driven} presentation.
His asymptotic methods are not so easy to trace: certainly, some early contact with A. Odlyzko played a role, and what was then often nick-named 
\emph{\`a la Odlyzko} later became \emph{singularity analysis of generating functions} \cite{FlOd}. One of his favorite methods was the Mellin transform that, according to his own words, he picked up from Rainer Kemp \cite{FNebel}, but under the name \emph{Gamma function method}. However, he quickly read parts of Doetsch's book \cite{Doetsch} on integral transforms. He also had a strong liking of Ramanujan and loved Hardy's book \cite{hardy-lectures}.

Later in life, Flajolet became more of a systematic theory-builder and concentrated more on (complex) analysis than on combinatorics. Much of his legacy can be found in his
book (with Sedgewick) \emph{Analytic Combinatorics} \cite{FS}. 

This historic paper  tries to emphasize his earlier years and some of his combinatorial results. From a personal perspective, the years 1980--1995 were the ones where I learnt most from him, and that is the period that will be described here in some way.

Several topics, like continued fractions, generation of combinatorial objects, Boltzmann samplers, urn models etc. could not be discussed.

\section{Register function (Horton-Strahler numbers)}

Flajolet's research on the register function of binary trees as published in
  the final article~\cite{FRV} was his first work
in the analysis of algorithms. This problem is nontrivial 
but manageable, and one can learn  a lot while working on it. 
As in most of Flajolet's work, the first part is of a combinatorial nature, while the second deals more with the asymptotic evaluation.

The problem deals with binary trees (counted by Catalan numbers) and the parameter \textsf{reg},
which associates to each binary tree (which is used to code an arithmetic expression,
with data in the leaves and operators in the internal nodes) the minimal number of extra registers that is needed to evaluate the tree. The optimal strategy is to evaluate difficult subtrees first, and use one register to keep its value, which does not hurt, if the other subtree requires less registers. If both subtrees are equally difficult, then one more register is used, compared to the requirements of the subtrees.

There is a recursive description of this function: $\textsf{reg}(\square)=0$, and if tree $t$ has subtrees $t_1$ and $t_2$, then
\begin{equation*}
	\textsf{reg}(t)=
	\begin{cases}
		\max\{\textsf{reg}(t_1),\textsf{reg}(t_2)\}&\text{ if } \textsf{reg}(t_1)\ne\textsf{reg}(t_2),\\
		1+\textsf{reg}(t_1)&\text{ otherwise}.
	\end{cases}
\end{equation*}

The register function is also known as Horton-Strahler numbers in the study of the complexity of river networks. The original papers are~\cite{Horton45, Strahler52};
since then many papers have been written on the subject.
The recursive description attaches numbers to the nodes, starting with 0's at the leaves and
then going up; the number appearing at the root is the register function of the tree.
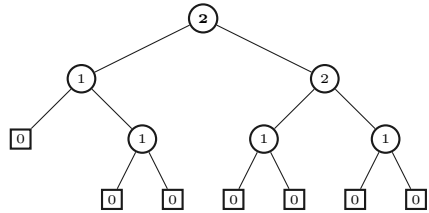
\begin{center}\tiny
	\begin{tikzpicture}
		[scale=0.4,inner sep=0.7mm,
		s1/.style={circle,draw=black!90,thick},
		s2/.style={rectangle,draw=black!90,thick}]
		\node(a) at ( 0,8) [s1] [text=black]{$\boldsymbol{2}$};
		\node(b) at ( -4,6) [s1] [text=black]{$1$};
		\node(c) at ( 4,6) [s1] [text=black]{$2$};
		\node(d) at ( -6,4) [s2] [text=black]{$0$};
		\node(e) at ( -2,4) [s1] [text=black]{$1$};
		\node(f) at ( 2,4) [s1] [text=black]{$1$};
		\node(g) at ( 6,4) [s1] [text=black]{$1$};
		\node(h) at ( -3,2) [s2] [text=black]{$0$};
		\node(i) at ( -1,2) [s2] [text=black]{$0$};
		\node(j) at ( 1,2) [s2] [text=black]{$0$};
		\node(k) at ( 3,2) [s2] [text=black]{$0$};
		\node(l) at ( 5,2) [s2] [text=black]{$0$};
		\node(m) at ( 7,2) [s2] [text=black]{$0$};
		\path [draw,-,black!90] (a) -- (b) node{};%
		\path [draw,-,black!90] (a) -- (c) node{};%
		\path [draw,-,black!90] (b) -- (d) node{};%
		\path [draw,-,black!90] (b) -- (e) node{};%
		\path [draw,-,black!90] (c) -- (f) node{};%
		\path [draw,-,black!90] (c) -- (g) node{};%
		\path [draw,-,black!90] (e) -- (h) node{};%
		\path [draw,-,black!90] (e) -- (i) node{};%
		\path [draw,-,black!90] (f) -- (j) node{};%
		\path [draw,-,black!90] (f) -- (k) node{};%
		\path [draw,-,black!90] (g) -- (l) node{};%
		\path [draw,-,black!90] (g) -- (m) node{};%
	\end{tikzpicture}
\end{center}
%\end{figure}

Let $\mathscr{R}_{p}$ denote the family of
trees with register function $=p$, then one gets immediately from the recursive 
definition:

\begin{center}\small
	\begin{tikzpicture}
		[inner sep=1.3mm,
		s1/.style={circle=10pt,draw=black!90,thick},
		s2/.style={rectangle,draw=black!50,thick},scale=0.5]
		
		\node at ( -5,0) { $\mathscr{R}_p$};
		
		\node at (-4,0) { $=$};
		\node(a) at (-2,1)[s1]{};
		\node(b) at (-3,-1){ $\mathscr{R}_{p-1}$};
		\node(c) at (-1,-1){ $\mathscr{R}_{p-1}$};
		\path [draw,-,black!90] (a) -- (b) node{};%
		\path [draw,-,black!90] (a) -- (c) node{};%
		\node at (0.7,0) {$+$};
		%\node at (1.5,0) {$2$};
		\node(d) at (3,1)[s1]{};
		\node(e) at (2,-1){ $\mathscr{R}_{p}$};
		\node(f) at (4,-1.2){ $\sum\limits_{j<p}\mathscr{R}_{j} $};
		\path [draw,-,black!90] (d) -- (e) node{};%
		\path [draw,-,black!90] (d) -- (f) node{};%
		\node at (5+0.7,0) {$+$};
		%5\node at (5+1.5,0) {$2$};
		\node(dd) at (5.5+3,1)[s1]{};
		\node(ee) at (5.5+2,-1.2){ $\sum\limits_{j<p}\mathscr{R}_{j}$};
		\node(ff) at (5.5+4,-1){ $\mathscr{R}_{p}$};
		\path [draw,-,black!90] (dd) -- (ee) node{};%
		\path [draw,-,black!90] (dd) -- (ff) node{};%
	\end{tikzpicture}
\end{center}

Later in life, Flajolet would write such \emph{symbolic equations} in a linearized form.
In terms of generating functions, these equations read as
\begin{equation*}
	R_p(z)=zR_{p-1}^2(z)+2zR_p(z)\sum_{j<p}R_j(z)=\frac{1-u^2}{u}\frac{u^{2^p}}{1-u^{2^{p+1}}}
	\quad\text{with}\quad z=\frac{u}{(1+u)^2};
\end{equation*}
the variable $z$ is used to mark the size (i.~e., the number of internal nodes) of the binary tree.
The substitution into the $u$-world originated in an important paper by de Bruijn, Knuth, and Rice~\cite{BrKnRi}.

The explicit formula for $R_p(z)$ can be proved by induction. This was of course not the way it was found.
Nowadays, with computer algebra, from a list of the first few values of $R_p(z)$, expressed in the variable $u$, one would
quickly guess the formula.

Explicit forms for 
$[z^n]R_p(z)$, the number of binary trees of size $n$ with register function $=p$, are available via contour integration (always around a small circle around the origin):
\begin{align*}
[z^n]R_p(z)&=\frac1{2\pi i}\oint\frac{dz}{z^{n+1}}\frac{1-u^2}{u}\frac{u^{2^p}}{1-u^{2^{p+1}}}
=[u^n](1+u)^{2n+2}\frac{(1-u)}{(1+u)^3}\frac{1-u^2}{u}\frac{u^{2^p}}{1-u^{2^{p+1}}}
\\&=\sum_{k\ge0}\bigg[\binom{2n}{n+1-(2k+1)2^p}-2\binom{2n}{n-(2k+1)2^p}+\binom{2n}{n-1-(2k+1)2^p}\bigg].
\end{align*}
From this, the average value of the register function of a random binary tree with $n$ nodes can be computed:
\begin{equation} \label{v2}
	[z^n]\sum_{p\ge1}pR_p(z)=\sum_{k\ge1}v_2(k)\bigg[\binom{2n}{n+1-k}-2\binom{2n}{n-k}+\binom{2n}{n-1-k}\bigg]
\end{equation}
with $v_2(k)$ being the number of trailing zeroes in the binary representation of $k$.
This is an interesting number theoretic quantity described by the recursive rules $v_2(2k)=1+v_2(k)$, $v_2(2k+1)=0$, $v_2(1)=0$.
Let $S_2(n)$ be the sum-of-digits function of $n$ in the binary number system. Then $S_2(n)-S_2(n-1)=1-v_2(n)$, and by telescoping
$S_2(n)=n-\sum_{k\le n}v_2(k)$. A further summation leads to a result by Delange~\cite{Delange}, published just before Flajolet's work on the register function was begun:
$\sum_{m<n}S_2(m)=\frac n2\log_2n+nF(\log_2n)$, where $F(x)$ is a continuous periodic function of period 1 with known Fourier coefficients. Using Abel's partial summation twice
on $[z^n]R_p(z)$, one gets eventually:

The average number of registers to evaluate a binary tree with $n$ 
nodes is asymptotically given by
\begin{equation*}
	\log_4n+D(\log_4n)+o(1)
\end{equation*}
with
\begin{equation*}
	D(x)=\sum_{k\in\mathbb{Z}}d_ke^{2\pi i kx}
\quad\text{and}\quad
	d_0=\frac12-\frac{\gamma}{2\log2}-\frac1{\log2}+\log_2\pi,\
	d_k=\frac1{\log2}\zeta(\chi_k)\Gamma(\tfrac{\chi_k}{2})(\chi_k-1),\quad k\ne0,
\end{equation*}
with $\chi_k=\frac{2\pi i k}{\log2}$. The classical Gamma- and zeta-functions appear here. 

At the same time, Rainer Kemp \cite{Kemp} proved equivalent results using different methods.

\section{{Digital search trees and Tries}}

Digital search trees are constructed from (sufficiently long) sequences of 0 and 1 (appearing with equal probability in the simplest model).
The order in which the items are inserted is important. In the following example, the items are inserted in the order $A,B,C,\dots$.
The most important parameter is the \emph{insertion depth}, i. e., the average distance of a node to the root. In \cite{Vol3}, there was an open problem that
Flajolet and Sedgewick \cite{FSdigital} solved: What is the (average) number of internal endnotes, if a digital search tree of $n$ is considered. In the example, 
$C,F,I,H$ are such endnodes, and the count is 4. We will describe the combinatorial part of the solution. The asymptotic part is based on Rice's method and
is perhaps the method of choice when it comes to analyze Digital search trees and their closely related cousins, called Tries.

\begin{minipage}[t] {0.38\linewidth}

%\begin{figure}[h!]
%	\begin{center}
		\begin{tikzpicture}[
		s1/.style={
			circle,
			very thick,
			draw=black,
			top color=white,
			bottom color=yellow!00!black!00
		},
		s2/.style={
			circle,
			draw=white!50!black!80,
			top color=white,
			bottom color=yellow!00!black!00
		},scale = 0.75]
		\path node at (1,0) {$A:1001$};%
		\path node at (1,-.5) {$B:0110$};%
		\path node at (1,-1) {$C:0000$};%
		\path node at (1,-1.5) {$D:1111$};%
		\path node at (1,-2) {$E:0100$};%
		\path node at (1,-2.5) {$F:0101$};%
		\path node at (1,-3) {$G:1101$};%
		\path node at (1,-3.5) {$H:1110$};%
		\path node at (1,-4) {$I\;:1100$};%
		\footnotesize %
		\path node(a) at (7,0) [s2, text centered]{$A$};%
		\path node(b) at (5,-1.5) [s2, text centered]{$B$};%
		\path node(c) at (4,-3) [s2, text centered]{$C$};%
		\path node(d) at (9,-1.5) [s2, text centered]{$D$};%
		\path node(e) at (6,-3) [s2, text centered]{$E$};%
		\path node(f) at (5,-4.5) [s2, text centered]{$F$};%
		\path node(g) at (10,-3) [s2, text centered]{$G$};%
		\path node(h) at (11,-4.5) [s2, text centered]{$H$};%
		\path node(i) at (9,-4.5) [s2, text centered]{$I$};%
		\path [draw,-,black!90] (a) -- (b) node[above,pos=.6,black]{{\tiny $0$}};%
		\path [draw,-,black!90] (b) -- (c) node[above,pos=.6,black]{{\tiny $0$}};%
		\path [draw,-,black!90] (a) -- (d) node[above,pos=.6,black]{{\tiny $1$}};%
		\path [draw,-,black!90] (b) -- (e) node[above,pos=.6,black]{{\tiny $1$}};%
		\path [draw,-,black!90] (e) -- (f) node[above,pos=.6,black]{{\tiny $0$}};%
		\path [draw,-,black!90] (d) -- (g) node[above,pos=.6,black]{{\tiny $1$}};%
		\path [draw,-,black!90] (g) -- (h) node[above,pos=.6,black]{{\tiny $1$}};
		\path [draw,-,black!90] (g) -- (i) node[above,pos=.6,black]{{\tiny $0$}};%
		\end{tikzpicture}
%	\end{center}

%\end{figure}
\end{minipage}

\hskip 9 cm
\begin{minipage}[t]{0.42\linewidth}
	\vskip -4.0cm  
	\noindent This is a \emph{digital search tree} with $9$ nodes. A \emph{trie} is very similar, except that the data is located in the leaves, not the internal nodes. Parameters like \emph{size} and \emph{path-length} have interesting connections to $q$-analysis.
	
\end{minipage}
%-------------------------------------------------------------------

Let 
\begin{equation*}
F_n(z)=\sum_{k\ge0}[\text{probability that a random digital search tree of size $n$ has $k$ endnodes}]z^k.
\end{equation*}
It is easy to see that
\begin{equation*}
F_{n+1}(z)=\sum_{k=0}^n2^{-n}\binom nk  F_k(z)F_{n-k}(z), \quad n\ge1, \quad F_0(z)=1,\quad F_1(z)=z,
\end{equation*}
since $2^{-n}\binom nk$ is the probability that $k$ items go to the left and $n-k$ go the right. The desired average is then obtained via
$\ell_n=F_n'(1)$, following the recursion
\begin{equation*}
\ell_{n+1}=2^{1-n}\sum_{k=0}^n\binom nk  \ell_k, \quad n\ge1, \quad \ell_0=0,\quad \ell_1=1.
\end{equation*}
This type of recursions is usually solved using exponential generating functions, since it translates into a functional-differential equation:
\begin{equation*}
L(z)=\sum_{n\ge0}\ell_n\frac{z^n}{n!}\qquad\text{and}\qquad L'(z)=1+2L\Bigl(\frac z2\Bigr)e^{z/2}.
\end{equation*}
This becomes better using the \emph{Poisson transformed function}:
\begin{equation*}
	\widehat L(z)=\sum_{n\ge0}\widehat\ell_n\frac{z^n}{n!}=e^{-z}L(z)\qquad\text{and}\qquad \widehat L'(z)=e^{-z}+2\widehat L\Bigl(\frac z2\Bigr).
\end{equation*}
Reading off the coefficients of $z^n/n!$ on both sides, we are led to the recursion
\begin{equation*}
\widehat\ell_{n+1}=(-1)^n-(1-2^{1-n})\widehat\ell_n,\quad n\ge0,\quad \widehat\ell_0=0,\qquad \text{which can be solved by iteration. }
\end{equation*}
Setting
\begin{equation*}
Q_n=\prod_{1\le \ell\le n}\Bigl(1-\frac1{2^k}\Bigr), \quad\text{then}\quad \widehat{\ell}_n=(-1)^{n+1}Q_{n-2}\Bigl(\frac1{Q_0}+\cdots+\frac1{Q_{n-2}}\Bigr),
\quad n\ge2,\quad \widehat{\ell}_0=0,\ \widehat{\ell}_1=1.
\end{equation*}
Since $L(z)=e^z\widehat{L}(z)$, comparing  coefficients eventually leads to the formula
\begin{equation*}
l_n=n-\sum_{k\ge2}\binom nk(-1)^kR_{k-2}, \quad\text{with} \quad R_n=Q_n\Bigl(\frac1{Q_0}+\cdots+\frac1{Q_{n}}\Bigr), \quad\text{for}\ n\ge2.
\end{equation*}
For the asymptotic evaluation of this with Rice's method, it is essential to have a definition of $R_n$, where $n$ can also be a complex number; currently the definition only makes sense for nonnegative integers.

Although not much will be done in this historic paper with Rice's method, at least the key formula should be mentioned:
\begin{equation*}
\sum_{k=1}^n\binom nk(-1)^kf(k)=\frac1{2\pi i}\int_{\mathscr{C}}\frac{(-1)^{n-1}n!}{z(z-1)\dots (z-n)},
\end{equation*}
where the positively oriented curve $\mathscr{C}$ encircles the poles $1,\dots,n$ and no others. The lower index 1 in the summation and contour could be replaced by
any number; in our application, 2 works. For that, the sequence $f(k)$ must be extended to an analytic function $f(z)$. This is often easy, but sometimes not.
What one usually does (if it is feasible) is to change the contour and to compute some extra residues.

It is interesting to observe that $Q_n$ really belongs to $q$-analysis. Using the notation $(z;q)_n=(1-z)(1-qz)\dots(1-q^{n-1}z)$, one could write $Q_n=(\frac12;\frac12)_n$.
Using ideas from $q$-analysis, the finite product can be written as a quotient of two infinite products, and that leads already to the definition of $Q_n$ for complex $n$:
\begin{equation*}
Q_n=\frac{Q_{\infty}}{Q(2^{-n})}, \quad\text{with} \quad Q(x):=\prod_{k\ge1}\Bigl(1-\frac x{2^k}\Bigr)\quad\text{and} \quad Q_\infty=Q(1)=0.288788\dots.
\end{equation*}
For $R_n/Q_n$, this is not so easy. A first idea is to write the finite sum as a difference of two infinite series; the problem is that $\frac1{Q_0}+\frac1{Q_1}+\cdots$ does not converge. In order to overcome this, one must pull out something so that the strategy works for the (tiny) rest. Skipping a few details, what works is
\begin{equation*}
R_n=n+1-\alpha +R^*_n, \quad\text{with} \quad \alpha:=\sum_{k\ge1}\frac1{2^k-1} \quad\text{and} \quad
R^*(z)=\sum_{j\ge0}\frac{(z+1+j-\alpha)2^{-z-1-j}}{(1-2^{-z-1})(1-2^{-z-2})\dots(1-2^{-z-1-j})},
\end{equation*}
then $R^*_n=R^*(n)$. In this series, $z$ can be a complex number; there is however a nicer version: 
\begin{equation*}
R^*(z)=\frac1{Q_\infty}\sum_{j\ge1}\frac{(-1)^{j-1}2^{-\binom j2}}{Q_{j-1}}\frac{z+j}{2^{z+j}-1},\quad\text{which can be obtained by partial fraction decomposition }
\end{equation*}
and a partition identity of Euler
\begin{equation*}
\sum_{n\ge0}\frac{t^n}{(q;q)_n}=\frac{1}{(1-t)(1-qt)(1-q^2t)\dots}=\frac1{(t;q)_\infty}. \quad\text{The special values $t=q=\frac12$ evaluate the series}
\end{equation*}
\begin{equation*}
\sum_{\ell\ge0}\frac{2^{-\ell}}{Q_\ell}=\frac1{Q_\infty}, \quad\text{and, after differentiation,}\quad
\sum_{\ell\ge0}\frac{\ell2^{-\ell}}{Q_\ell}=\frac{\alpha}{Q_\infty}.
\end{equation*}
The rest of the method is of a more analytic nature and shall not be discussed here. Eventually, it turns out that the average number of endnodes is (apart from small fluctuations)
$\ell_n\approx n\cdot 0.372048\dots$; the constant can be analytically expressed as $\alpha+1-R^*(-1)$. Recall that in the term for $j=1$ in $R^*(-1)$  a limit must be taken.

\section{Approximate Counting}

In \cite{BIT} the following scenario is investigated. There is a counter $C$, initially set to $1$. Assume that the counter's value is currently  $k$ and a random element arrives, then the counter advances to value $k+1$ with probability $2^{-k}$.

After $n$ random increments, the counter has a value between $1$ and $n$, but is typically around $\log_2n$. So the counter serves as an approximate count, with significantly less bits.

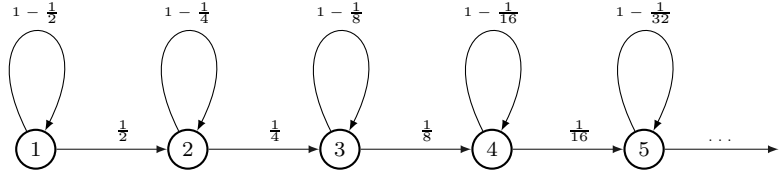
\begin{center}
	\begin{tikzpicture}[
		s1/.style={
			circle,
			thick,
			draw=black,
			top color=white,
			bottom color=yellow!00!black!00
		},
		s2/.style={
			circle,
			draw=white!50!black!80,
			top color=white,
			bottom color=yellow!50!black!30
		}]
		\begin{scope}[>=latex]
			\footnotesize %
			\path node(a) at (0,0) [s1, text centered]{$1$};%
			\path node(b) at (2,0) [s1, text centered]{$2$};%
			\path node(c) at (4,0) [s1, text centered]{$3$};%
			\path node(d) at (6,0) [s1, text centered]{$4$};%
			\path node(e) at (8,0) [s1, text centered]{$5$};%
			
			\path [draw,->,black!90] (a) -- (b) node[above,pos=.6,black]{{\tiny $\frac12$}};%
			\path [draw,->,black!90] (b) -- (c) node[above,pos=.6,black]{{\tiny $\frac14$}};%
			\path [draw,->,black!90] (c) -- (d) node[above,pos=.6,black]{{\tiny $\frac18$}};%
			\path [draw,->,black!90] (d) -- (e) node[above,pos=.6,black]{{\tiny $\frac1{16}$}};%
			%\draw [->](-0.23,0.23) .. controls (-1,2) and (1,2) .. (0.23,0.23);
			\path[draw,->](a) .. controls (-1,2) and (1,2) .. (a)node[above,pos=0.5,black]{{\tiny $1-\frac12$}};
			\path[draw,->](b) .. controls (1,2) and (3,2) .. (b)node[above,pos=0.5,black]{{\tiny $1-\frac14$}};
			\path[draw,->](c) .. controls (3,2) and (5,2) .. (c)node[above,pos=0.5,black]{{\tiny $1-\frac18$}};
			\path[draw,->](d) .. controls (5,2) and (7,2) .. (d)node[above,pos=0.5,black]{{\tiny $1-\frac1{16}$}};
			\path[draw,->](e) .. controls (7,2) and (9,2) .. (e)node[above,pos=0.5,black]{{\tiny $1-\frac1{32}$}};
			\path [draw,->,black!90] (e) -- (10-0.23,0) node[above,pos=.5,black]{{\tiny $\dots$}};%
		\end{scope}
	\end{tikzpicture}
\end{center}
The graph described this very well, and also allows to find
the generating function
	\begin{align*}
	H_\ell(z)&=\sum_{n\ge0}z^n[\text{Probability to reach state $\ell$ after $n$ random steps}]=\frac{z^{\ell-1}q^{\binom \ell 2}}{\prod\limits_{i=1}^\ell\bigl(1-(1-q^i)z\bigr)}
\end{align*}
easily. Partial fraction decomposition and reading off coefficients leads to:
\begin{equation*}
	p_{n,\ell}:=[z^n]H_\ell(z)=\sum_{i=0}^{\ell-1}\frac{(-1)^i2^{-\binom i2}}{Q_iQ_{\ell-1-i}}\Big(1-2^{-(\ell -i)}\Big)^n\quad\text{with}
	\quad Q_i:=\Big(1-\frac1{2}\Big)\Big(1-\frac1{4}\Big)\dots\Big(1-\frac1{2^i}\Big);
\end{equation*}
 the products $Q_i$ appeared already in the section on digital search trees.
Expanding the binomial in $p_{n,\ell}$ according to the binomial theorem, and using Euler's partition identities
\begin{equation*}
	\prod_{m\ge0}\big(1+{x}{q^m}\big)=\sum_{j\ge0}\frac{x^jq^{\binom j2}}{(q;q)_j},\qquad
\prod_{m\ge0}\big(1-{x}{q^m}\big)^{-1}=\sum_{j\ge0}\frac{x^j}{(q;q)_j}
\end{equation*}
one can find the average value $C_n$ of the counter after $n$ random increments to be
\begin{equation*}
	C_n=1-\sum_{k=1}^n (-1)^{k}\binom nk2^{-k}Q_{k-1}
	=1-\frac1{2\pi i}\int_{\mathscr C}\frac{(-1)^{n-1}n!}{z(z-1)\dots(z-n)}
	2^{-z}Q_{z-1}f(z) dz.
\end{equation*}
The contour $\mathscr C$ encircles the poles $1,2,\dots,n$ and no others.
The appearance of such contour integrals is typical for Rice's method. Eventually, using methods from complex analysis one finds that $C_n\approx\log_2n$. The study of the lower order (fluctuating) terms is very instructive, but shall not be presented here.

There are other ways to compute the probabilities $p(n,k)$. In \cite{LP08}, this is done by iteration:
\begin{equation*}
\text{set}\quad F(z,u):=\sum_{n\ge0}\sum_{k\ge1}p(n,k)z^nu^k, \quad \text{then one derives ($q=\frac12$ here)}
\end{equation*}
\begin{equation*}
F(z,u)=\frac{u}{1-z}+\frac{z(u-1)}{1-z}F(z,qu)=\sum_{j\ge0}\frac{(-1)^jz^j(q;u)_juq^j}{(1-z)^{j+1}}.
\end{equation*}

\section{Probabilistic Counting}

Assume that $n$ data are given and, via a \emph{hash-function}, each data is mapped to a non-negative integer, following the geometric distribution with parameter $q=\frac12$.
We might think about $n$ persons, and to each person a ball is associated, with a number printed on it, and the probabilities are $\frac12,\frac14,\frac18,\dots$.
Then each person throws her ball into an urn with the corresponding number. We expect that about $\frac n2$ balls fall into urn $0$, $\frac n4$ balls fall into urn $1$,
and so on. The goal of the paper \cite{FMartin} is to use these urns to produce a count of the number of data (persons) without actually keeping track of them.
Since we only distinguish which urn is empty or which urn is non-empty, it does not matter whether a person appears at the counter more than once, since the number associated to her is always the same. 

The parameter $R_n$ that is used is the index of the first empty urn, or the longest sequence of non-empty urns (starting with urn 0). Flajolet computed $q_{n,k}$,
the probability that $R_n\ge k$, using an inclusion-exclusion argument. This is somewhat surprising since he told me that he does not like inclusion/exclusion too much, unless it corresponds to generating functions where a factor $(1-z)^k$ generates the alternating signs automatically.

 So we need to work out the probability that urns $0,1,\dots,k-1$ are non-empty.
This is given by 1 minus the probability that a particular urn is empty, plus the probability that two particular urns are empty, and so on. We get
\begin{equation*}
q_{n,k}=1-\sum_{0\le i <k}\Bigl(1-\frac1{2^i}\Bigr)^n+\sum_{0\le i<j <k}\Bigl(1-\frac1{2^i}-\frac1{2^j}\Bigr)^n	
-\sum_{0\le i<j<\ell <k}\Bigl(1-\frac1{2^i}-\frac1{2^j}-\frac1{2^\ell}\Bigr)^n+\cdots.
\end{equation*}
Using different indices, we can rewrite this:
\begin{equation*}
	q_{n,k}=\sum_{t=0}^{k-1}(-1)^t\sum_{0\le i_0<\dots<i_{t-1}<k}\Bigl(1-\frac{2^{i_0}+\cdots 2^{i_{t-1}}}{2^k}\Bigr)^n.
\end{equation*}
We note that $2^{i_0}+\cdots+ 2^{i_{t-1}}$ is just a (unique) non-negative integer, written in binary, and the sign depends on the number of bits in the binary representation.
Denote by $\nu(k)$ the number of digits 1 in the binary representation of $k$ (it was denoted by $S_2(k)$ in a previous section), then we have the beautiful formula
\begin{equation*}
	q_{n,k}=\sum_{j=0}^{2^k}(-1)^{\nu(j)}\Bigl(1-\frac{j}{2^k}\Bigr)^n.
	\quad \text{The sequence $(-1)^{\nu(j)}$ is known as the Prouhet-Thue-Morse-sequence.}
\end{equation*}
The pattern of signs goes like $+--+-++--++-+--+ \dots$. The average value of the parameter $R_n$ is exactly given by $\overline{R_n}=\sum_{k\ge 1}q_{n,k}$. Using methods from real analysis that will not be repeated here, Flajolet approximates:
\begin{equation*}
q_{n,k}\approx \psi\Bigl(\frac{n}{k}\Bigr)\quad\text{with}\quad \psi(x)=\prod_{j\ge0}(1-e^{-x2^j})=\sum_{j\ge0}(-1)^{\nu(j)}e^{-jx}.
\end{equation*}
The last expansion is classical and easier to read when setting $q:=e^{-x}$. Continuing his Mellin transform analysis, Flajolet is lead to study an interesting function, which
resembles the Riemann $\zeta$-function:
\begin{equation*}
N(s)=\sum_{j \ge 1 }\frac{(-1)^{\nu(j)}}{j^s};\quad\text{$N(s)$ should be read `nu-of-s', but  there is no special character available.}
\end{equation*} 
One major question is the analytic continuation of $N(s)$; as it stands, it only converges for $\Re(s)>1$. The trick is to group for consecutive terms together, according to the pattern $+--+$; then there is enough cancellation so that one can go to the left with the variable $s$:
\begin{equation*}
N(s)=-1-2^{-s}+3^{-s}+\sum_{j\ge0}\frac{(-1)^{\nu(j)}}{(4j)^s}\biggl[1-\Bigl(1+\frac{1}{4j}\Bigr)^{-s}-\Bigl(1+\frac{2}{4j}\Bigr)^{-s}+
\Bigl(1+\frac{3}{4j}\Bigr)^{-s}\biggr].
\end{equation*}
Eventually one gets that the average value of the parameter $R_n$ is asymptotically given by
\begin{equation*}
\overline{R_n}\approx \log_2(\varphi\cdot n)+\text{tiny fluctuation,\ \ and the constant $\varphi$ involves } \prod_{p\ge1}\biggl[\frac{(4p+1)(4p+2)}{(4p)(4p+3)}\biggr]^{(-1)^{\nu(p)}}.
\end{equation*}
Such products are not uncommon in the literature \cite{AllSha}, but not too much is known about the Flajolet-Martin constant \cite{Finch-book}.

\section{Adding a new slice}

In \cite{FP}, the sequence $1,1,2,3,5,9,16,28,\dots$ (A002572 in \cite{Sloane}) is investigated used the adding-a-new-slice technique. This technique was probably known to some people under no name or a different name, but Flajolet's nickname for it is very descriptive. I vaguely remember that even Cayley was interested in this sequence, but I lost the reference.

The definition of the \emph{level number sequences } is
\begin{equation*}
\text{C1. }n_1=1,\qquad \text{C2. For all $j$ such that } 1< j\le k:\ n_j\le 2n_{j-1} ,\qquad \text{C3. } n=n_1+\cdots+n_k.
\end{equation*}
The number $k$ is called the height and the number of representations of $n$ created in this way is denoted by $H_n$. A generating function $F(q,u)$ is formed where $q$ keeps track of the current total value and $u$ of the last entry. One writes $F(q,u)=F_1(q,u)+F_2(q,u)+\cdots$, and the index refers to the height of the sequence. The process to go from 
$F_{k}$ to $F_{k+1}$ is called adding a new slice. It is driven by the substitution
\begin{equation*}
u^j\longrightarrow uq+(uq)^2 +\cdots +(uq)^{2j}=\frac{uq}{1-uq}\Big(1-(uq)^{2j}\Big), \quad\text{or}\quad F_{k+1}(q,u)= \frac{uq}{1-uq}\Big(F_k(q,1)-F_k(q,u^2q^2)\Big).
\end{equation*}
Summing on $k$ leads to
\begin{equation*}
 F(q,u)=uq+ \frac{uq}{1-uq}\Big(F(q,1)-F(q,u^2q^2)\Big),\quad\text{which can be solved by iteration:}
\end{equation*}
\begin{equation*}
 F(q,u)=uq+ \frac{uq}{1-uq}F(q,1)-\frac{uq}{1-uq}\Big[(uq)^2 +\frac{(uq)^2}{1-(uq)^2}F(q,1) -\frac{(uq)^2}{1-(uq)^4}\Big[(uq)^4+\frac{(uq)^4}{1-(uq)^4}F(q,1) + \cdots.
\end{equation*}
Plugging in $u=1$ and collecting leads to the explicit representation
\begin{equation*}
F(q,1)=\dfrac{\sum\limits_{j\ge1}(-1)^{j+1}\frac{q^{2^{j+1}-j-2}}{(1-q)(1-q^3)(1-q^7)\dots(1-q^{2^{j-1}-1})}}
{1-\sum\limits_{j\ge1}(-1)^{j+1}\frac{q^{2^{j+1}-j-2}}{(1-q)(1-q^3)(1-q^7)\dots(1-q^{2^{j}-1})}}.
\end{equation*}
Recall that $F(q,1)$ is the generating function of the numbers $H_n$, since setting $u=1$ means that the height $k$ is no longer of relevance. Asymptotics are not difficult from here (but not discussed in detail), since it can be shown that there is a dominant simple pole at $0.5573678719$. Eventually this leads to
\begin{equation*}
H_n=0.25450\cdot 1.794147^n+O(1.43^n); \quad\text{the exponential growth constant is the reciprocal of the   simple pole.}
\end{equation*}

The author had a chance to use the adding-a-new-slice technique several times after the event of the level numbers paper.

\section{Harmonic numbers and Euler sums}

It is known that the quantity $\displaystyle{S_n(m)=\sum_{k=1}^n}\binom nk\frac{(-1)^k}{k^m}$ for a negative integer $m$ is basically a subset Stirling number (`of the second kind'). Somebody asked Flajolet what happens for \emph{positive} integers $m$? Here is the answer, expressed in terms of \emph{harmonic numbers} 
$H_n^{(j)}=\sum_{1\le k\le n}\frac1{k^j}$: We start from a (cleverly chosen) quantity $T$ and use partial fraction decomposition:
\begin{equation*}
T=\frac{n!}{z(z-1)(z-2)\dots(z-n)}\frac1{z^m}=\sum_{k=1}^n\binom nk(-1)^{n-k}\frac1{k^m}\frac1{z-k}+\frac{\lambda}{z^{m+1}}+\cdots+\frac{\mu}{z}.
\end{equation*}
Now we consider $z\cdot T$ and let $z\to\infty$ with the result
\begin{equation*}
\sum_{k=1}^n\binom nk(-1)^{n-k}\frac1{k^m}+{\mu}=0; \quad \mu=[z^{-1}]T=[z^{m}]\frac{n!}{(z-1)\dots(z-n)}
=(-1)^n[z^m]\frac1{(1-z)(1-\frac z2)\dots(1-\frac zn)}.
\end{equation*}
The computation continues:
\begin{equation*}
\mu=-\sum_{k=1}^n\binom nk(-1)^{n-k}\frac1{k^m}=(-1)^n[z^m]\exp\log\frac1{(1-z)(1-\frac z2)\dots(1-\frac zn)}
=(-1)^n[z^m]\exp\bigg\{\sum_{k=1}^n\sum_{j\ge1}\frac1{jk^j}\bigg\},
\end{equation*}
or
\begin{align*}
\sum_{k=1}^n\binom nk(-1)^{k-1}\frac1{k^m}
	&= [z^m]\exp\bigg\{\sum_{j\ge1}\frac{H_n^{(j)}}{j}\bigg\}\\
	&=[z^m]\biggl(1+H_nz+\Big(\frac{H_n^2}{2}+\frac{H_n^{(2)}}{2}\Big)z^2+
\Big(\frac{H_n^3}{6}+\frac{H_nH_n^{(2)}}{2}+\frac{H_n^{(3)}}{3}\Big)z^3+\cdots\biggr),
\end{align*}
which provides the first three evaluations:
\begin{equation*}
H_n,\quad \frac{H_n^2}{2}+\frac{H_n^{(2)}}{2}, \quad \frac{H_n^3}{6}+\frac{H_nH_n^{(2)}}{2}+\frac{H_n^{(3)}}{3};
\end{equation*}
a general formula can be written using (partial) Bell polynomials. The paper \cite{FSfinite} contains more such alternating summations, where
rational functions no longer work, and instead of exact evaluations one gets asymptotic evaluations using contour integrations and Rice's integrals. ---

Harmonic numbers appear again, for instance in \cite{FSalvy} when studying \emph{Euler sums}
\begin{equation*}
S_{p,q}=\sum_{n=1}^\infty\frac{H_n^{(p)}}{n^q},\quad\text{like}\quad S_{1,2}=2\zeta(3),\ S_{1,3}=\frac54\zeta(4),\ S_{2,4}=\zeta(3)^2-\frac13\zeta(6),\ \&c.
\end{equation*}
The technique of the highly cited paper \cite{FSalvy} is \emph{residue calculus}, writing the sum in question in two different ways as a sum of residues using a suitable
\emph{kernel} function.

\section{Flajolet and Ramanujan}

A famous unproven assertion by Ramanujan is:
\begin{equation*}
\frac12e^n=1+\frac{n!}{1!}+\frac{n!}{2!}+\cdots+\frac{n!}{n!}\theta\quad\text{with}\quad \theta=\frac13+\frac{4}{135(n+k)}, \quad \frac8{45}<k<\frac2{21}. 	
\end{equation*}
This was proved in \cite{FGKT}: There are auxiliary functions
\begin{equation*}
Q(n)=1+\frac{n-1}{n}+\frac{(n-1)(n-2)}{n^2}+\cdots\quad\text{and}\quad R(n)=1+\frac{n}{n+1}+\frac{n^2}{(n+1)(n+2)}+\cdots;
\end{equation*}
it is not difficult to prove that  $Q(n)+R(n)=n!e^n/n^n$; $Q(n) $ is called Ramanujan's $Q$-function. It turns out that Ramanujan's assertion can be rephrased as
\begin{equation*}
R(n)-Q(n)=\frac23+\frac{8}{135(n+k)}, \quad\text{with }\quad \frac8{45}<k<\frac2{21}.
\end{equation*}
Instead of real analysis which was used in earlier attempts, complex analysis will be used, namely a famous implicit function defined by $y=ze^y$; it either goes by the name
Lambert's $W$ function or tree function. The generating function has the explicit expansion
\begin{equation*}
y=\sum_{n\ge1}n^{n-1}\frac{z^n}{n!};\quad  \text{recall that $n^{n-1}$ is the number of rooted labelled trees with $n$ nodes.}
\end{equation*} 
 One finds a generating function
\begin{equation*}
\sum_{n\ge1}Q(n)n^{n-1}\frac{z^n}{n!}\quad\text{and from it the asymptotic expansion}\quad Q(n)=\sqrt{\frac{\pi n}{2}}-\frac13+\cdots.
\end{equation*}
However, for Ramanujan's assertion, \emph{explicit} error bounds are requested. This is done in the above mentioned paper using contour integration and careful bounds and computer algebra. ---

Flajolet was also asked about his personal hit-parade of Ramanujan's results \cite{Notices}. I am not sure what he answered, but I know for sure that he was very fond
of \emph{Ramanujan's master theorem} \cite{master}, which can be stated as
\begin{equation*}
\int_0^\infty x^{s-1}\big(\lambda(0)-x\lambda(1)+x^2\lambda(2)-\cdots\big)dx=\frac{\pi }{\sin(\pi s)}\lambda(-s)\quad\text{for}\quad
0<\Re(s)<1,
\end{equation*}
under some technical conditions on  $\lambda$. This, among other things, has inspired Flajolet to come up with his Poisson-Newton-Mellin-Rice cycle \cite{cycle}.

\section{Digits}

The existence of equation (\ref{v2}) motivated Flajolet (later joined by coauthor Ramshaw) \cite{FRamshaw} to look for something similar in a problem solved by Sedgewick around the same time \cite{Batcher}. Without describing the odd-even merge algorithm, it turns out that
\begin{equation*}
\sum_{k\ge1}\theta(k)\binom{2n}{n-k}\quad\text{with}\quad \theta(k)=\begin{cases}
	\ 1&\text{ if } k=2^i(4j+1)\\
	-1&\text{ if } k=2^i(4j+3)\\
\end{cases}
\end{equation*}
is the quantity to be studied. Note that each positive integer $k$ has a unique representation in this form.
The quantity $\theta(k)$ plays the role of $v_2(k)$. Double summation does not lead to the binary representation,
but to something similar, namely the Gray code representation of integers. This representation might be described via
\begin{equation*}
(n)_{\text{GR}}=\dots a_2a_1a_0 \quad\text{with}\quad a_k=\biggl\lfloor\frac{n}{2^{k+2}}+\frac 34\biggr\rfloor
-\biggl\lfloor\frac{n}{2^{k+2}}+\frac 14\biggr\rfloor.
\end{equation*}
The Gray code representation has the pleasant property that from $n-1$ to $n$ exactly one bit changes. This means that for the sum-of-digits 
function $S_{\text{GR}}(n)-S_{\text{GR}}(n-1)=\pm1$, and, more specifically, 
$S_{\text{GR}}(n)-S_{\text{GR}}(n-1)=\theta(n)$. The telescoping summation leads to
$S_{\text{GR}}(n)=\sum_{m\le n}\theta(m)$. A further summation leads to a beautiful explicit result of the Delange type \cite{FRamshaw}.

 In this way, Sedgewick's Mellin type analysis was
replaced by a more elementary approach. The typical periodic fluctuations appear, as well as  the Dirichlet series $\zeta(s,\tfrac34)-\zeta(s,\tfrac14)$. ---

Later, as part of his series about applications of the Mellin transformation in Combinatorics, Flajolet (plus four coauthors) \cite{FGKPT} proposes an elegant 
high level approach to such questions related to all kinds of question involving various digits representations: 

\textbf{The Mellin-Perron formula:} Let $c>0$ lie in the half-plane of absolute convergence of the Dirichlet series
$\sum_k\lambda k^{-s}$. Then for any $m\ge1$, we have
\begin{equation*}
\frac1{m!}\sum_{1\le k <n}\lambda_k\Big(1-\frac kn\Big)^m=
\frac1{2\pi i}\int_{c-i\infty}^{c+i\infty}
\Big(\sum_{k\ge1}\frac{\lambda_k}{k^s}\Big)n^s
\frac{ds}{s(s+1)\dots(s+m)}.
\end{equation*}
For $m=0$,
\begin{equation*}
\sum_{1\le k <n}\lambda_k+\frac{\lambda_n}{2}
=\frac1{2\pi i}\int_{c-i\infty}^{c+i\infty}
\Big(\sum_{k\ge1}\frac{\lambda_k}{k^s}\Big)n^s
\frac{ds}{s}.
\end{equation*}
The instance $m=0$ is usually called Perron's summation formula. For problems with the sum-of-digits function, we usually need $m=1$, which is fortunate, since there are less issues about convergence. In most cases, the sum in question can be evaluated by residues. The result is either asymptotic or even exact, as in Delange's case, or in the instance of the Gray code representation. Without going into details, the formula
\begin{equation*}
0=\int_{-\frac14-i\infty}^{-\frac14+i\infty}
\zeta(s)n^s\frac{ds}{s(s+1)}
\end{equation*}
is responsible for the exact formulae. The interested reader is referred to \cite{FGKPT}.

\bibliographystyle{plain}
%\bibliography{historic}

\end{document}